\documentclass[12pt]{amsart}
\usepackage{%amsmath, 
amssymb, %amsthm, 
epsfig}
\theoremstyle{plain}
 
  \newtheorem{theorem}{Theorem}[subsection]%[section]%

\theoremstyle{remark}

\numberwithin{equation}{section}
\newcommand{\resq}{\overline{\mathcal{U}}_q(sl_2)}
\newcommand{\sresq}{\widehat{\mathcal{U}}_q(sl_2)}

\newcommand{\bw}{{{w}}}
\newcommand{\End}{{\operatorname{End}}}
\begin{document}

\markboth{J. Murakami, K. Nagatomo}
{Logarithmic knot invariants}

%\catchline{}{}{}{}{}

\title{LOGARITHMIC KNOT INVARIANTS ARISING FROM RESTRICTED QUANTUM GROUPS}

\author{JUN MURAKAMI}

\address{Department of Mathematics
\\
Faculty of Science and Engineering
\\
 Waseda University
 \\
 3-4-1 Ohkubo, Shinjyuku-ku, Tokyo 169-8555, JAPAN. 
\\
E-mail:\ murakami@waseda.jp}

\author{KIYOKAZU NAGATOMO}

\address{Department of Pure and Applied Mathematics
\\
Graduate School of Information Science and Technology
\\
 Osaka University
 \\ Toyonaka, Osaka 560-0043, JAPAN. 
 }

\maketitle

%\baselineskip18pt
%\begin{center}
%\begin{Large}\textbf{
%Logarithmic knot invariants arising from restricted quantum groups}
%\end{Large}
%\vskip0.5cm
%Jun Murakami\footnote{Department of Mathematics, School of Science and Engineering, Waseda University, 3-4-1 Ohkubo, Shinjyuku-ku, Tokyo 169-8555, JAPAN. 
%E-mail:\ murakami@waseda.jp}
% and Kiyokazu Nagatomo\footnote{Department of Pure and Applied Mathematics, Graduate School of Information Science and Technology, Osaka University, Toyonaka, Osaka 560-0043, JAPAN.}
 %, Koji Ohnuki
%\end{center}
%\bigskip
\begin{abstract}
We construct knot invariants from the radical part of projective modules of the restricted quantum group 
$\resq$ at $q = \exp(\pi \sqrt{-1}/p)$,
and  we also show a relation between these invariants and the colored Alexander invariants.  
These projective modules are related to logarithmic conformal field theories.  
\end{abstract}
\section{Introduction}
Various knot invaraiants are constructed from the quantum $R$-matrix of the quantum group.  
However, most of them are constructed from semisimple algebras.  
Our concern in this note is constructing knot invariant arising from {\it non-semisimple} representations.   
We focus on the restricted quantum group $\resq$ and construct knot invariant which is understood as a derivative of the colored Alexander invariant \cite{ADO}, \cite{M}.  
\par
Let $Z$ and $J$ be the center and the Jacobson radical of  $\resq$ respectively.  
Then $Z$ is a direct sum of $Z^{(s)}$ and $Z^{(r)}$, where $Z^{(s)}$ is the subalgebra of $Z$ generated by the primitive idempotents 
and $Z^{(r)} = Z \cap J$.  
Let $K$ be a knot in $S^3$.  
By using the idea of the universal invariant in \cite{O}, 
we can associate an element $z_K \in Z$ with $K$.   
Then $z_K$ is expressed as
$$
z_K = z_K^{(s)} + z_K^{(r)}   \qquad 
(z_K^{(s)} \in Z^{(s)}, \ z_K^{(r)} \in Z^{(r)}).
$$
The first term $z_K^{(s)}$  is a linear combination of the primitive idempotents and
the coefficient of each idempotent corresponds to the colored Jones invariant.   
In this paper,  we study about the knot invariants coming from $z_K^{(r)}$.  
The space $Z^{(r)} = Z \cap J$ has a natural basis corresponding to the indecomposable modules, and the coefficients of $z_K^{(r)}$ with respect to this basis are also knot invariants.  
\par
In the construction of $z_K$, we assume that $K$ is a single component knot.  
For a multi-component link, the construction in this paper does not work well and we need another idea to extend $z_K^{(r)}$ for link case.   
\par
A three-manifold invariant is constructed from the `integral' of $\resq$ \cite{Ku}, \cite{H}, \cite{KR}, \cite{O2}, which is defined for a finite dimensional Hopf algebra.    
Meanwhile, an action of $SL(2, {\mathbb Z})$ on the center $Z$ of $\resq$ is given in \cite{Ker} and \cite{F1}.  
By using  our invariants constructed here, 
these two theories can be combined as the usual topological quantum field theory, e.g.  \cite{T}, \cite{BK},  related to ${\mathcal U}_q(sl_2)$.  
The detail will be given elsewhere.  
\par
We review the definition of  $\resq$ and its representations in Section 2.  
The construction of $z_K^{(r)}$ is given  in Section 3.
In Sections 4 and 5, we show some property of the invariants coming from  $z_K$, 
especially the relation to the colored Alexander invariant in 
\cite{M}.  
\section{Restricted quantum groups}
\subsection{Definition}
Let $p \geq 2$ be a positive integer and $q = \exp(\pi \sqrt{-1}/p)$.  
The semi-restricted quantum group  $\sresq$
is the quotient of the usual quantum group $\mathcal{U}_q(sl_2)$ defined by the following generators and relations
as an algebra.  
$$
\begin{aligned}
&\sresq = \left<
K, K^{-1}, E, F \mid
K \, K^{-1} = K^{-1}\, K = 1, 
\right.\\
&\left.
K \, E \, K^{-1} = q^2 \, E, \ 
K \, F \, K^{-1} = q^{-2} \, F, \ 
E \, F - F \, E = \dfrac{K - K^{-1}}{q-q^{-1}}, 
\right. \\
&\left.
E^p = F^p = 0
\right>.  
\end{aligned}
$$
The restricted quantum group $\resq$ is obtained from $\sresq$ by inquiring one more relation $K^{2p}=1$.  
The coproduct, counit and antipode of $\sresq$ and $\resq$ are defined as follows.  
$$
\begin{aligned}
\Delta(K) = K \otimes K,
\quad
\Delta(E) 
&=
E \otimes K + 1 \otimes E, \quad
\Delta(F)
=
F \otimes 1 + K^{-1} \otimes F, 
\\
&\epsilon(E)= \epsilon(F) =0, \quad
\epsilon(K) = 1, 
\\
S(E) &= -E \, K^{-1}, \quad
S(F) = -K \, F, \quad
S(K) = K^{-1}.  \qquad
\end{aligned}
$$
\subsection{$R$-matrix}
By introducing a symbol $k$ such that $k^2 = K$, 
We can define an $R$-matrix $R$ of $\resq$ satisfying
$$
R \, \Delta(X) \, R^{-1} = \overline{\Delta}(X),
$$
where
$
\overline{\Delta}(X) = \sum_i z_i^{} \otimes y_i^{}
$
if $\Delta(X) = \sum_i y_i^{} \otimes z_i^{}$.
The explicit form of $R$ is given as follows.  
\begin{equation}
R 
=
q^{\frac{1}{2}H \times H} 
\sum_{n=0}^{p-1}
\dfrac{(q-q^{-1})^n}{[n]!} \, 
q^{\frac{n(n-1)}{2}} \,
(E^n \otimes F^n), 
\label{e:u}
\end{equation}
where $[n] = \frac{q^n - q^{-n}}{q-q^{-1}}$ and $[n]! = [n][n-1]\cdots[1]$, and $H$ is an element such that $q^H = K$.  
%\subsection{Ribbon element}
%The restricted quantum group $\resq$ does not have $R$-matrix in itself because $k \notin \resq$, but has the ribbon element $v$ and is given actually in \cite{F1} as follows. 
%$$
%v
%=
%\dfrac{1-i}{2\sqrt{p}}
%\sum_{m=0}^{p-1}\sum_{j=0}^{2p-1}
%\dfrac{(q-q^{-1})^m}{[m]!}
%q^{-m/2+mj+(j+p+1)^2/2} F^mE^mK^j
%$$ 
\subsection{Irreducible mudules}
The irreducible modules ${\mathcal X}^\alpha(s)$ of $\resq$ are labeled by $\alpha=\pm1$ and $1\leq s \leq p$, and is spanned by weight vectors 
 $\left|s, n\right>^{\pm}$, $0 \leq n \leq s-1$ with the action of 
 $\resq$ given by
$$
\begin{aligned}
K \left|s, n\right>^{\pm} &= 
\pm q^{s-1-2n}\left|s, n\right>^{\pm},
\\
E  \left|s, n\right>^{\pm} &= 
\pm [n][s-n]\left|s, n-1\right>^{\pm},
\\
F \left|s, n\right>^{\pm} &= 
\left|s, n+1\right>^{\pm},
\end{aligned}
$$
where $\left|s, s\right>^{\pm} = 0$.  
\subsection{Projective modules}
Projective modules ${\mathcal P}^\pm(s)$  of $\resq$ which are fundamental to investigate the structure of the center of $\resq$ 
are labeled by $1 \leq s \leq p-1$.  
Note that a $\resq$ module is also naturally a $\sresq$ module. 
\par
Let $s$ be any integer $1 \leq s \leq p-1$.
The projective module ${\mathcal P}^+(s)$ has the basis
$$
\{x_k^{(+, s)}, y_k^{(+, s)}\}_{0 \leq k \leq p-s-1}
\cup
\{a_n^{(+,s)}, b_n^{(+, s)} \}_{0\leq n\leq s-1},
$$
and the action of $\resq$ is given by
$$
\begin{aligned}
K \, x_k^{(+, s)}
&=
-q^{p-s-1-2k} \, x_k^{(+, s)}, 
\qquad
K \, y_k^{(+, s)}
=
-q^{p-s-1-2k} \, y_k^{(+, s)}, 
\qquad
0 \leq k \leq p-s-1,
\\
K \, a_n^{(+, s)}
&=
q^{s-1-2n} \, a_n^{(+, s)}, 
\qquad
K \, b_n^{(+, s)}
=
q^{s-1-2n} \, b_n^{(+, s)}, 
\qquad
0 \leq n \leq s-1,
\\
E\, x_k^{(+, s)} 
&=
-[k][p-s-k] \, x_{k-1}^{(+,s)},
\qquad
0 \leq k \leq p-s-1, \quad
(\text{with} \ x_{-1}^{(+, s)} = 0),
\end{aligned}
$$
$$
\begin{aligned}
E\, y_k^{(+, s)}
&=
\begin{cases}
-[k][p-s-k]\,y_{k-1}^{(+,s)}, & 1 \leq k \leq p-s-1, \\
a_{s-1}^{(+,s)}, & k=0,
\end{cases}
\\
E\, a_n^{(+, s)}
&=
[n][s-n]\,a_{n-1}^{(+,s)}, 
\quad 0 \leq n \leq s-1, \quad
(\text{with}\quad a_{-1}^{(+,s)} = 0), 
\\
E\, b_n^{(+, s)}
&=
\begin{cases}
[n][s-n]\,b_{n-1}^{(+,s)}+a_{n-1}^{(+,s)}, 
& 1 \leq n \leq s-1, \\
x_{p-s-1}^{(+,s)}, & n=0,
\end{cases}
\end{aligned}
$$
and
$$
\begin{aligned}
F\, x_k^{(+, s)}
&=
\begin{cases}
x_{k+1}^{(+,s)}, & 0 \leq k \leq p-s-2, \\
a_{0}^{(+,s)}, & k=p-s-1,
\end{cases}
\\
F\, y_k^{(+, s)}
&=
y_{k+1}^{(+,s)}, \quad 0\leq k \leq p-s-1 ,
 \quad
(\text{with} \ y_{p-s}^{(+, s)} = 0),
\\
F\, a_n^{(+, s)}
&=
a_{n+1}^{(+,s)}, 
\quad 0 \leq n \leq s-1, \quad
(\text{with}\quad a_{s}^{(+,s)} = 0), 
\\
F\, b_n^{(+, s)}
&=
\begin{cases}
b_{n+1}^{(+,s)}, &
0 \leq n \leq s-2, \\
y_0^{(+, s)}, & n=s-1.   
\end{cases}
\end{aligned}
$$
\par
Let $s$ be an integer $1 \leq s \leq p-1$.
The projective module ${\mathcal P}^-(p-s)$ has the basis
$$
\{x_k^{(-, s)}, y_k^{(-, s)}\}_{0 \leq k \leq p-s-1}
\cup
\{a_n^{(-,s)}, b_n^{(-, s)} \}_{0\leq n\leq s-1},
$$
and the action of $\resq$ is given by
$$
\begin{aligned}
K \, x_k^{(-, s)}
&=
-q^{p-s-1-2k} \, x_k^{(+, s)}, 
\qquad
K \, y_k^{(+, s)}
=
-q^{p-s-1-2k} \, y_k^{(-, s)}, 
\qquad
0 \leq k \leq p-s-1,
\\
K \, a_n^{(-, s)}
&=
q^{s-1-2n} \, a_n^{(-, s)}, 
\qquad
K \, b_n^{(-, s)}
=
q^{s-1-2n} \, b_n^{(-, s)}, 
\qquad
0 \leq n \leq s-1,
\\
E\, x_k^{(-, s)} 
&=
-[k][p-s-k] \, x_{k-1}^{(-,s)},
\qquad
0 \leq k \leq p-s-1, \quad
(\text{with} \ x_{-1}^{(-, s)} = 0)
\end{aligned}
$$
$$
\begin{aligned}
E\, y_k^{(-, s)}
&=
\begin{cases}
-[k][p-s-k]\,y_{k-1}^{(-,s)} + x_{k-1}^{(-,s)}, & 1 \leq k \leq p-s-1, \\
a_{s-1}^{(-,s)}, & k=0,
\end{cases}
\\
E\, a_n^{(-, s)}
&=
[n][s-n]\,a_{n-1}^{(-,s)}, 
\quad 0 \leq n \leq s-1, \quad
(\text{with}\quad a_{-1}^{(-,s)} = 0), 
\\
E\, b_n^{(-, s)}
&=
\begin{cases}
[n][s-n]\,b_{n-1}^{(-,s)}, 
& 1 \leq n \leq s-1, \\
x_{p-s-1}^{(-,s)}, & n=0,
\end{cases}
\end{aligned}
$$
and
$$
\begin{aligned}
F\, x_k^{(-, s)}
&=
x_{k+1}^{(-,s)}, \quad 0 \leq k \leq p-s-1,
\quad
(\text{with} \ x_{p-s}^{(-, s)} = 0), \\
F\, y_k^{(-, s)}
&=
\begin{cases}
y_{k+1}^{(-,s)}, & 0\leq k \leq p-s-2 , \\
b_0^{(-,s)}, & k = p-s-1,
\end{cases}
\\
F\, a_n^{(-, s)}
&=
\begin{cases}
a_{n+1}^{(-,s)}, &
 0 \leq n \leq s-2,
\\
x_0^{(-,s)}, & n=s-1,
\end{cases}
\\
F\, b_n^{(-, s)}
&=
b_{n+1}^{(-,s)}, \quad
0 \leq n \leq s-1,
\quad
(\text{with}\quad b_{s}^{(-,s)} = 0).
\end{aligned}
$$
Note that the diagonal part is a direct sum of irreducible modules.  
\subsection{Center}
The dimension of the center ${\mathcal Z}$  of $\resq$ is $3p-1$.
% and
%is spanned by the centers of irreducible representations and the %radical parts of the projective representations ${\mathcal P}^\pm(s)$.  
The basis of ${\mathcal Z}$ is given by the canonical central elements in \cite{F1} as follows:
Two special primitive idempotents $e_0$ and $e_p$, other primitive idempotents $e_s$, $(1 \leq s \leq p-1)$, and $2\,(p-1)$ elements $w_s^\pm$ $(1 \leq s \leq p-1)$ corresponding to the radical part.  
These basis satisfy the following relations. 
\begin{equation}
\begin{aligned}
e_s \, e_{s'} &= \delta_{s, s'} \, e_s,   &
s, \  s' = 0, \ldots, p, 
\\
e_s \, w_{s'}^\pm &= \delta_{s, s'} \, w_s^\pm, &
\qquad
0 \leq x \leq p, \ 1 \leq s' \leq p-1,
\\
w_s^\pm \, w_{s'}^\pm &= w_s^\pm \, w_{s'}^\mp = 0, &
\quad 1 \leq s, s' \leq p-1.  
\end{aligned}
\label{e:centerrelation}
\end{equation}
\section{Logarithmic invariants of knots}
\subsection{Knots and  (1,1)-tangles}
\begin{figure}[htb]
\begin{center}
\epsfig{file=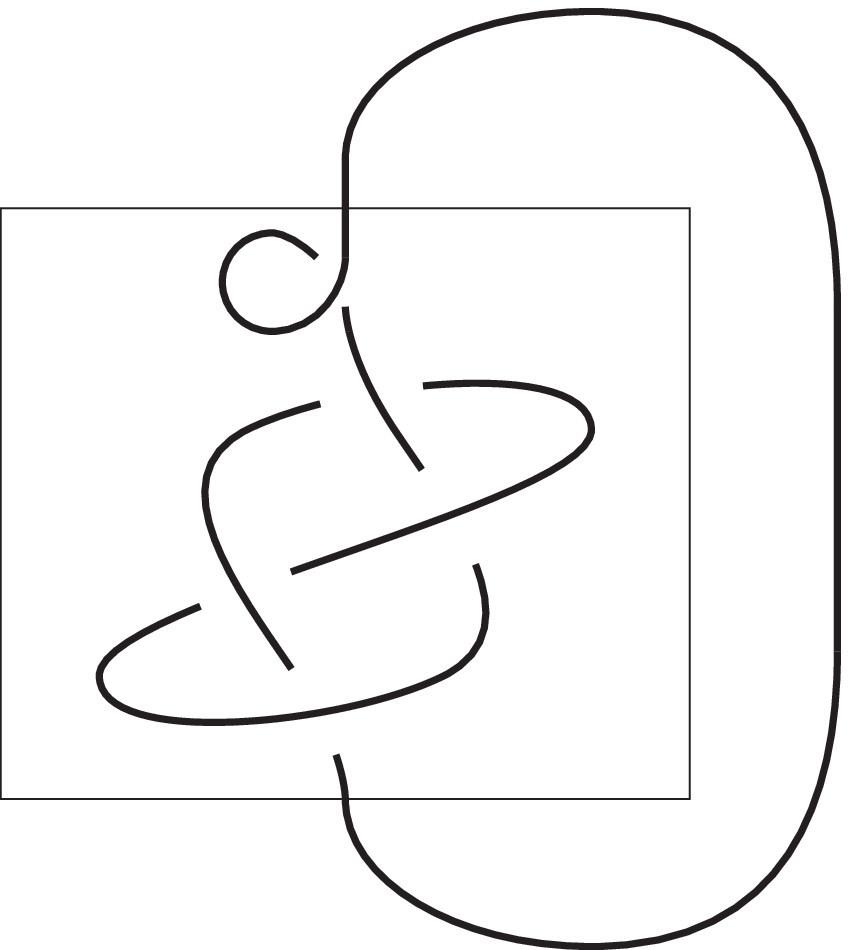, scale=0.4}
\qquad\qquad
\raisebox{4mm}{
\epsfig{file=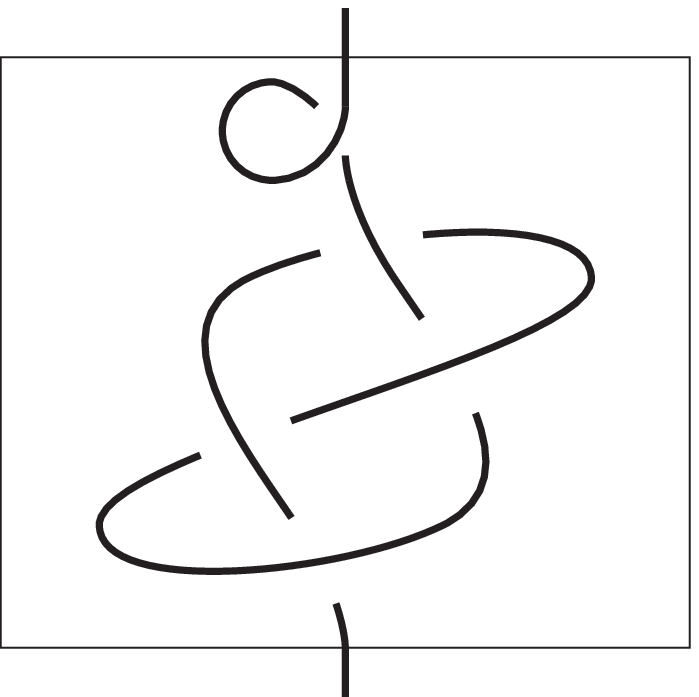, scale=0.4}
}
\\
$K$  \qquad\qquad\qquad\qquad\qquad\qquad $T_K$
\\
Figure eight knot $K$ with framing 1 
and its tangle $T_K$. 
\end{center}
\caption{Cosure of a framed tangle.}
\label{figure:tangle}
\end{figure}

In this paper, knots and tangles are oriented and framed.  
For a connected (1,1)-tangle $T$, let $K_T$ be the knot obtained by  joining the two open ends as in Figure \ref{figure:tangle}.  
For two tangles $T$ and $T'$, it is known that $K_T$ and $K_{T'}$ are isotopic as framed knots if and only if $T$ and $T'$ are isotopic as framed tangles.  
So, in the following, we sometimes mix up invariants of connected (1,1)-tangles and  invariants of knots. 
\subsection{Framed braid}
\begin{figure}[htb]
\begin{center}
$
\begin{matrix}
\sigma_i \hfill& {\sigma_i}^{-1} \hfill&
\tau_i \hfill& {\tau_i}^{-1}\hfill
\\
\quad
\raisebox{9.5mm}{$\cdots$}
\ %
\epsfig{file=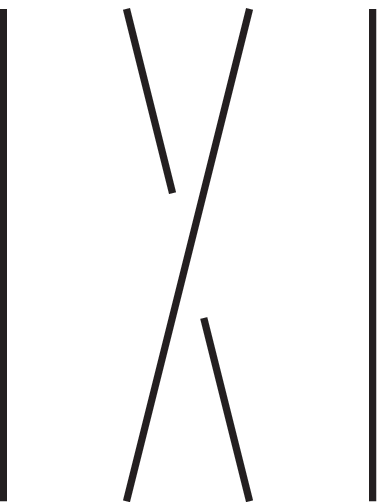, scale=0.4}
\ %
\raisebox{9.5mm}{$\cdots$}
\quad
&
\quad
\raisebox{9.5mm}{$\cdots$}
\ %
\epsfig{file=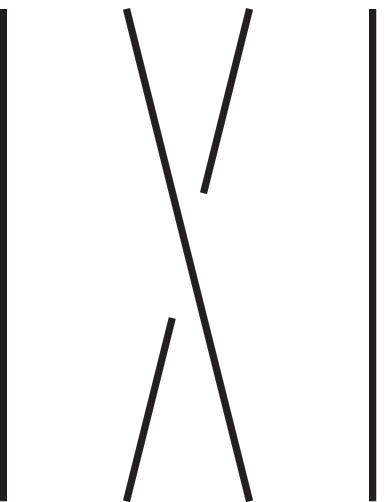, scale=0.4}
\ %
\raisebox{9.5mm}{$\cdots$}
\quad
&
\quad
\raisebox{9.5mm}{$\cdots$}
\ %
\epsfig{file=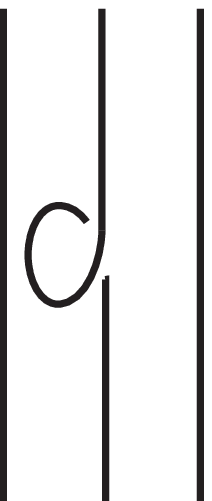, scale=0.4}
\ %
\raisebox{9.5mm}{$\cdots$}
\quad
&
\quad
\raisebox{9.5mm}{$\cdots$}
\ %
\epsfig{file=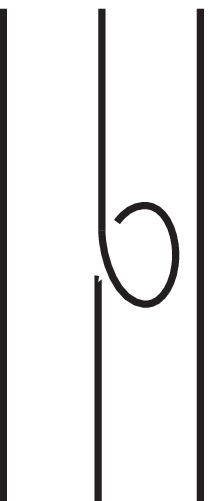, scale=0.4}
\ %
\raisebox{9.5mm}{$\cdots$}
\quad
\\
\quad i \ i+1 &\quad i \ i+1 & i & i
\end{matrix}
$
\end{center}
\caption{Generators of the framed braid group $FB_n$.}
\label{figure:generator}
\end{figure}

Framed braid group on $n$ strings $FB_n$ is defined by the following generators and relations.  
$$
\begin{aligned}
FB_n
&=
\left<
\sigma_1^{}, \sigma_2^{}, \dots, \sigma_{n-1}, 
\tau_1^{}, \tau_2^{}, \dots, \tau_n
\mid
\right.
\\& \qquad\qquad
\sigma_i^{}\, \sigma_{i+1}^{} \, \sigma_i^{} = 
\sigma_{i+1}^{} \, \sigma_i^{} \, \sigma_{i+1}, \quad
(i = 1, 2, \ldots, n-2),
\\& \qquad\qquad
\sigma_i^{} \, \sigma_j^{}
=
\sigma_j^{} \, \sigma_i^{}, \quad (|i-j| > 1),
\\& \qquad\qquad
\tau_i^{\pm1} \, \sigma_i^{} = \sigma_i^{} \, \tau_{i+1}^{\pm1}
,\quad
%\quad (i = 1, 2, \ldots, n-1)
%\\& \qquad\qquad
\tau_{i+1}^{\pm1} \, \sigma_i^{} = \sigma_i^{} \, \tau_{i}^{\pm1},
\quad (i = 1, 2, \ldots, n-1),
&\\& \qquad\qquad
\sigma_i^{} \, \tau_j^{}
=
\tau_j^{} \, \sigma_i^{}
, \quad
% \quad (|i-j| > 1)
%\\& \qquad\qquad
\left.
\tau_i^{} \, \tau_j^{}
=
\tau_j^{} \, \tau_i^{} \quad (|i-j| > 1)
\right>.
\end{aligned}
$$
The generators $\sigma_i^{\pm1}$ correspond to the positive and negative crossings, and $\tau_i^{\pm1}$ represent the blackboard framing corresponding to the twist as in Figure \ref{figure:generator}.  
Let $S_n$ be the symmetric group of $n$ letters $\{1$, 2, $\ldots$, $n\}$,
and $\pi$ be the group homomorphism from $FB_n$ to $S_n $ sending $\sigma_i^{}$ to the transposition $(i, i+1)$ and $\tau_i^{}$ to the identity for $i = 1$, 2, $\ldots$, $n-1$.  
\subsection{Alexander's and Markov's theorems}
Then we have the framed versions of Alexander's and Markov's theorem as follows.  
\begin{theorem}
Any framed link is isotopic to the closure $\hat b$ of some framed braid $b \in FB_n$.    
\end{theorem}
\begin{theorem}
Two framed braids $b_1^{}$, $b_2^{}$ have isotopic closures if and only if $b_1^{}$ can be transformed to $b_2^{}$ by a finite sequence of moves of the following two types.  
\begin{enumerate}
\item[(i)]
\quad
$\beta_1^{}\, \beta_2^{} \longleftrightarrow \beta_2^{} \, \beta_1^{}$
\quad
for $\beta_1^{}$, $\beta_2^{} \in FB_n$.  
\item[(ii)]
\quad
$\beta \, \tau_n^{\pm1} \longleftrightarrow i(\beta) \, \sigma_n^{\pm1}$
\quad
for $\beta\in FB_n \overset{i}\longrightarrow FB_{n+1}$.
\end{enumerate}
\end{theorem}
\subsection{Representation of $FB_n$ on $\bigotimes^n \resq$}
The universal $R$-matrix satisfies the Yang-Baxter equation
$$
R_{12} \, R_{13} \, R_{23}
=
R_{23} \, R_{13} \, R_{12}
:
\otimes^3\resq \longrightarrow \otimes^3\resq, 
$$
where $\otimes^3\resq = \resq \otimes \resq \otimes \resq$
and $R_{ij}$ acts on the $i$-th and $j$-th components of the tensor product.  
Let
$$
R = \sum_i r'_i \otimes r_i^{\prime\prime}
\ :\ 
\resq \otimes \resq \longrightarrow \resq \otimes \resq,  
$$
and
$$
v = \sum_i r_i^{\prime\prime} \, K^{N-1} \, r_i^{\prime}
\ :\ 
\resq \longrightarrow \resq.  
$$
Then $(\resq, R, v)$ is a ribbon Hopf algebra with  the ribbon element $v$.  
Therefore, we can define a homomorphism $\rho$ from $FB_n$ to $\End(\otimes^n \resq)$ by
$$
\begin{aligned}
\rho(\sigma_i^{})(x_1^{} \otimes \ldots \otimes x_i^{} \otimes x_{i+1}^{} \otimes \ldots \otimes x_n^{}) 
&= 
\sum_i
x_1^{} \otimes \ldots \otimes r_i^\prime \, x_{i+1}^{} \otimes 
r_i^{\prime\prime} \, x_{i}^{} \otimes \ldots \otimes x_n^{}, 
\\
\rho(\tau_i^{})(x_1^{} \otimes \ldots \otimes x_i^{} \otimes \ldots \otimes x_n^{})  
&= 
x_1^{} \otimes \ldots \otimes v \, x_i^{} \otimes \ldots \otimes x_n^{}.
\end{aligned}
$$
\subsection{Universal invariant}
Let $K$ be a knot and let $b_K^{}$ be a framed braid whose closure is equivalent to $K$.  
Then $\rho(b_K^{}) \in \End(\otimes^n \resq)$ is expressed as follows.
$$
\rho(b_K^{})
=
\sum_i a_{1, i}^{} \otimes a_{2, i}^{} \otimes \ldots \otimes a_{n, i}^{}.  
$$

Let $T$ be a (1,1)-tangle corresponding to $K$.  
The element $z_T \in \resq$ corresponding to $T$ is
defined by 
$$
z_T
=
a_{\pi(b_K)^{n-1}(1), i}\, K^{N-1}\,a_{\pi(b_K)^{n-2}(1), i}
\dots a_{\pi(b_K)^2(1), i}\, K^{N-1}\,a_{\pi(b_K)(1), i}\, K^{N-1}\,a_{1, i}.  
$$

% from $\tilde R$ in \cite{M}, which is an element of $\resq/I$ where $I$ is the ideal of $\resq$ generated by the elements of the form $\alpha\, \beta - \beta \, \alpha$.  
%Then $\phi(T)$ is in the center of $\resq$ and is an invariant of the framed knot $K$ 
%which is the closure of $T$, 
%\par
%?????
%????????
%\par
%???????
%??????
%\par
%and we also denote it by $\phi(K)$.   
%\par
The element $z_T \in \resq$ commutes with any elements of $\resq$ 
and is in the center ${\mathcal Z}$.  
% and is an invariant of the framed knot $K$ 
%which is the closure of $T$.  
Therefore, we have
\begin{equation}
z_T
=
\sum_{s=0}^p (a_s(T) \, {e}_s +
b_s^+(T)\, {w}_s^+ + b_s^-(T)\, \bw_s^-),  
\label{e:z}
\end{equation}
where 
$a_s(T)$, $b_s^\pm(T)$ are scalars and are invariants of the closure $K$ of $T$.  
Hence we can also denote them by $a_s(K)$, $b_s^\pm(K)$.  
\section{Properties of $a_s(K)$ and $b_s^\pm(K)$}
We show some property of   
$a_s(K)$ and $b_s^\pm(K)$.
\begin{theorem}
\begin{enumerate}
\item 
For the connected sum of two knots $K_1$ and $K_2$,   
$$
a_s(K_1 \# K_2) = a_s(K_1) \, a_s(K_2), 
\quad
b_s^\pm(K_1 \# K_2) = 
a_s(K_1) \, b_s^\pm(K_2) + b_s^\pm(K_1) \, a_s(K_2).  
$$
\item
For a knot $K$, the invariant $a_s(K)$ $(1 \leq s \leq p)$ is equal to the colored Jones invariant
$J_{s}(K)$ corresponding to the $s$-dimensional irreducible module ${\mathcal X}^+(s)$ normalized as 
$$
J_s({\text{unknot}}) = 1.
$$
\end{enumerate}
\end{theorem}
\begin{proof}
For two tangles $T_1$, $T_2$, let $T_1 \cdot T_2$ be the tangle obtained by joining $T_2$ below $T_1$.  
Then, for two knots $K_1$ and $K_2$,  
$T_{K_1} \cdot T_{K_2}$ is a tangle representing the connected sum $K_1 \# K_2$.  
Therefore, by using \eqref{e:centerrelation} and  \eqref{e:z}, 
we have
$$
\begin{aligned}
z_{K_1 \# K_2}
& =
z_{K_1} \, z_{K_2}
\\
&=
 \sum_{s=0}^p \left(a_s(K_1) \, a_s(K_2) \, {e}_s + 
 \right.
 \\&\!\!\!\!\!\!\!\!\!\!\!\!\!\!\!\!
 \left.
\left( 
b_s^+(K_1)\, a_s(K_2)+ a_s(K_1) \, b_s^+(K_2)
\right)\, {w}_s^+ 
+ 
\left( 
b_s^-(K_1)\, a_s(K_2)+ a_s(K_1) \, b_s^-(K_2)
\right)\, \bw_s^-
\right)   
\end{aligned}
$$  
and we obtain (1).  
\par
The center $e_s$ acts on ${\mathcal X}(s)$ as identity, and the other basis $e_t$ $(t \neq s)$, and $w_t^\pm$ acts on ${\mathcal X}(s)$ as zero.  
Hence $a_s(K)$ corresponds to the scalar representing the action of $z_K$ on ${\mathcal X}(s)$, which is equal to the colored Jones invariant and we get (2).  
\end{proof}
\section{Relation to the colored Alexander invariant}
\subsection{Relation}
Let $K$ be a framed knot and $T_K$ be the corresponding framed tangle.
Let $O_\lambda^p(T_K)$ be the scalar multiple of the colored Alexander invariant defined  in \cite{M}.  
Then we have the following.  
\begin{theorem}
The invariants $a_s(K)$, $b_s^+(K)$, $b_s^-(K)$ are given by the colored Alexander invariants as follows.  
$$
\begin{aligned}
a_s(K) &= O_{s-1}^p(T_K),
\qquad
0 \leq s \leq p, 
\\
b_s^+(T_K) 
&= 
- \dfrac{p\, \sin^2{\frac{\pi}{p}}}{\pi \,\sin{\frac{\pi s}{p}}} \, \left(\left.
\dfrac{d \, O_{\lambda}^p(T_K)}{d\lambda}\right|_{\lambda=2p-s-1}
-\left.
\dfrac{d \, O_{\lambda}^p(T_K)}{d\lambda}\right|_{\lambda=s-1}
\right),
\quad
1 \leq s \leq p-1,  
\\
b_s^-(K) 
&= 
\dfrac{p \, \sin^2{\frac{\pi}{p}}}{\pi \, \sin\frac{\pi \, s}{p}}  \, \left(\left.
\dfrac{d \, O_{\lambda}^p(T_K)}{d\lambda}\right|_{\lambda=s-1}
-\left.
\dfrac{d \, O_{\lambda}^p(T_K)}{d\lambda}\right|_{\lambda=-s-1}
\right),
\quad
1 \leq s \leq p-1.
\end{aligned}
$$
\end{theorem}
Before proving the above, we introduce some representations of 
$\sresq$.  
\subsection{Non-integral representations}
We introduce highest weight representations of $\sresq$ for 
non-integral weights  and 
obtain the projective modules ${\mathcal P}^{\pm}$ 
as a specialization of cetrain non-irreducible module.  
\par
First, we define the irreducible module for non-integer number $\lambda$ as follows.  
Let ${\mathcal X}(\lambda)$ be the $\sresq$ module spanned by weight vectors $v_n^\lambda$, $0 \leq n \leq p-1$.  
The action of $\sresq$ to ${\mathcal X}(\lambda)$ is given by 
$$
\begin{aligned}
K \, v_n^\lambda &= 
q^{\lambda-1-2n}\, v_n^\lambda,
\qquad
E \, v_n^\lambda = 
[n][\lambda-n]\, v_{n-1}^\lambda,
\qquad
F \, v_n^\lambda = 
\, v_{n+1}^\lambda,
\end{aligned}
$$
where $\, v_p^\lambda=0$.  
\par
Next, we define a non-irreducible module which is isomorphic to direct sum of two non-integral highest modules.  
Let $t$ be an integer with $1 \leq s \leq p$ and ${\mathcal Y}(\lambda, s)$ be the $\sresq$ module which is spanned by 
weight vectors $c_j^{(\lambda, s)}$ and $d_j^{(\lambda, s)}$ for $0 \leq j \leq p-1$.
The action of $\sresq$ is given by
$$
\begin{aligned}
K \, c_n^{(\lambda, s)}
&=
q^{\lambda-1-2n} \, c_n^{(\lambda, s)}, 
\qquad
K \, d_n^{(\lambda, s)}
=
q^{\lambda-1-2s-2n} \, d_n^{(\lambda, s)}, 
\qquad
0 \leq n \leq p-1,
\\
%\end{aligned}
%$$
%$$
%\begin{aligned}
E\, c_n^{(\lambda, s)}
&=
\begin{cases}
0, & \quad n = 0, \\
[n][\lambda-n]\,c_{n-1}^{(\lambda, s)} , 
&
\quad 1 \leq n \leq p-1 , 
\end{cases}
\\
E\, d_n^{(\lambda, s)}
&=
\begin{cases}
c_{s-1}^{(\lambda, s)}, & n = 0, \\
[n][\lambda-2s-n]\,d_{n-1}^{(\lambda, s)}
+ c_{n+s-1}^{(\lambda, s)},
& 1 \leq n \leq p-s, \\
[n][\lambda-2s-n]\,d_{n-1}^{(\lambda, s)}, 
& p-s+1 \leq n \leq p-1,
\end{cases}
\\
%\end{aligned}
%$$
%and
%$$
%\begin{aligned}
F\, c_n^{(\lambda, t)}
&=
\begin{cases}
c_{n+1}^{(\lambda, s)}, & 0 \leq n \leq p-2, \\
0, & n = p-1,
\end{cases}
\\
F\, d_n^{(\lambda, s)}
&=
\begin{cases}
d_{n+1}^{(\lambda, s)}, & 0 \leq n \leq p-2, \\
0, & n = p-1.
\end{cases}
\end{aligned}
$$
\subsection{Colored Alexander invariant}
Let $K$ be a framed knot, $T_K$ be the corresponding framed tangle,  and $\hat z_K$ be the corresponding central element in the semi-restricted quantum group $\sresq$ which is defined as $z_K$ by using the universal R-matrix of $\sresq$ given by \eqref{e:u}.  
Let $Z_K^{(\lambda, s)}$ be the representation matrix of $\hat z_K$ on ${\mathcal Y}(\lambda, s)$ with repsect to the above basis
$\{c_n^{(\lambda, s)}, d_n^{(\lambda, s)}; 0 \leq n \leq p-1\}$.  
Then the diagonal element corresponding to $c_n^{(\lambda, s)}$ and $d_n^{(\lambda, s)}$ $(0 \leq n \leq p-1)$ are 
equal to $O_{\lambda-1}^p(T_K)$ and  $O_{\lambda-1-2s}^p(T_K)$ respectively, where $O_\lambda^p(T_K)$ is given in \cite{M} as the scalar corresponding to the tangle $T_K$. 
Note that $O_\lambda^p(T_K)$ is a scalar multiple of the colored alexander invariant $\Phi_K^p(\lambda)$ and $O_\lambda^p(T_K)$ itself is also an invariant of $K$ if $K$ is a single component knot.    
\subsection{Proof of Theorem 4}
The matrix  $Z_K^{(\lambda, s)}$ has off-diagonal elements at
$(c_{n+s}, d_n)$ components for $0 \leq n \leq p-s$.  
Let $x$ be the $(c_{s}, d_0)$ component of 
 $Z_K^{(\lambda, s)}$.  
Then 
\begin{equation}
\begin{aligned}
Z_K^{(\lambda, s)}\,c_s^{(\lambda, s)}
&=
O_{\lambda-1}^p(T_K) \, c_s^{(\lambda, s)} ,
\\
Z_K^{(\lambda, s)}\,d_0^{(\lambda, s)}
&=
O_{\lambda-1-2s}^p(T_K) \, d_0^{(\lambda, s)} + 
x \, c_s^{(\lambda, s)}.  
\end{aligned}
\label{eq:1}
\end{equation}
Let
$$
h = c_{s}^{(\lambda, s)} - [s][\lambda-s] \, d_0^{(\lambda, s)}.  
$$
Then $E \, h=0$ and so $h$ is a highest weight vector of weight $\lambda-1 - 2s$.  
Therefore, on the one hand, 
\begin{equation*}
Z_K^{(\lambda, s)}\, h
=
O_{\lambda-1-2s}^p(T) \, h
=
O_{\lambda-1-2s}^p(T) \, c_s^{(\lambda, s)}
-
O_{\lambda-1-2s}^p(T) \,  [s][\lambda-s]\, d_0^{(\lambda, s)}
.  
\label{eq:2}
\end{equation*}
On the other hand, 
from \eqref{eq:1}, 
\begin{equation*}
Z_K^{(\lambda, s)}\, h
=
O_{\lambda-1}^p(T) \, c_0^{(\lambda, s)}
-
 [s][\lambda-s] \, x \, c_0^{(\lambda, s)}
-
O_{\lambda-1-2s}^p(T) \,  [s][\lambda-s]\, d_0^{(\lambda, s)}.
\label{eq:3}
\end{equation*}
Thus we have
\begin{equation}
O_{\lambda-2s-1}^p(T) = O_{\lambda-1}^p(T) - [s][\lambda-s] x, 
\label{eq:offdiagonal}
\end{equation}
and then
$$
x 
= 
\dfrac{O_{\lambda-1}^p(T)- O_{\lambda-1-2s}^p(T)}{[s][\lambda-s]}.  
$$
Hence, we get
\begin{equation}
\lim_{\lambda\to s+mp}\!\!  x
=
(-1)^m \, 
\dfrac{p\, \sin{\frac{\pi}{p}}}{\pi \, [s]} \, \left(\left.
\dfrac{d \, O_{\lambda-1}^p(T)}{d\lambda}\right|_{\lambda=s+mp}
\!\!\!\!
-\left.
\dfrac{d \, O_{\lambda-1-2s}^p(T)}{d\lambda}\right|_{\lambda=s+mp}
\right), 
\  %
m \in {\mathbf Z}.   
\label{eq:limit}
\end{equation}
The projective module ${\mathcal P}^+(s)$ is identical to ${\mathcal Y}(2p-s, p-s)$ by the correspondence of the basis 
$x_n^{(+, s)} \mapsto c_{n}^{(2p-s, p-s)}$, 
$a_n^{(+, s)} \mapsto c_{n+p-s}^{(2p-s, p-s)}$, 
$b_n^{(+, s)} \mapsto d_{n}^{(2p-s, p-s)}$, 
$y_n^{(+, s)} \mapsto d_{n+s}^{(2p-s, p-s)}$.   
Therefore, by substituting $2p-s$ to $\lambda$, $p-s$ to $s$, and $1$ to $m$ for \eqref{eq:limit}, we have
$$
a_s(T) = O_{2p-s-1}^p(T),
$$
and
$$
b_s^+(T) 
= 
-\dfrac{ p \, \sin^2{\frac{\pi}{p}}}{\pi \, \sin\frac{\pi \, s}{p}}  \,
\left(\left.
\dfrac{d \, O_{\lambda}^p(T)}{d\lambda}\right|_{\lambda=2p-s-1}
-\left.
\dfrac{d \, O_{\lambda}^p(T)}{d\lambda}\right|_{\lambda=s-1}
\right),
\quad
1 \leq s \leq p-1.   
$$
Similarly, 
 the projective module ${\mathcal P}^-(p-s)$ is identical to ${\mathcal Y}(s, s)$
by the correspondence of the basis 
$a_n^{(-, s)} \mapsto c_{n}^{(s, s)}$, 
$x_n^{(-, s)} \mapsto c_{n+s}^{(s, s)}$, 
$y_n^{(-, s)} \mapsto d_{n}^{(s, s)}$,    
$b_n^{(-, s)} \mapsto d_{n+p-s}^{(s, s)}$.  
 Hence we have
$$
b_{p-s}^-(T) 
= 
\dfrac{ p \, \sin^2{\frac{\pi}{p}}}{\pi \, \sin\frac{\pi \, s}{p}}  \,
\left( \left.
\dfrac{d \, O_{\lambda}^p(T)}{d\lambda}\right|_{\lambda=s-1}
-
\left.
\dfrac{d \, O_{\lambda}^p(T)}{d\lambda}\right|_{\lambda=-s-1}
\right),   
\quad
1 \leq s \leq p-1.   
$$
%\par
The formula \eqref{eq:offdiagonal} implies that
$$
O_{\lambda}^{mp+s-1}(T)
=
O_{\lambda}^{mp-s-1}(T).
$$
By putting $m=1$ and $s=p-s$, we have
$$
a_s^{}(T) 
=
O_{\lambda}^{2p-s-1}(T)
=
O_{\lambda}^{s-1}(T).
$$
\hfill
$\square$
\par\noindent


\begin{thebibliography}{abc}
\bibitem{ADO}
Y.  Akutsu, T. Deguchi and T. Ohtsuki, 
Invariants of colored links,
{\it J. Knot Theory Ramifications} {\bf 1} (1992) 161--184.  
\bibitem{BK}
 B. Bakalov and A. Kirillov Jr.,
 {\it Lectures on tensor categories and modular functors}, 
 University Lecture Series, {\bf 21} (American Mathematical Society, Providence, RI, 2001). 
\bibitem{F1}
B. L. Feigin, A. M.  Gainutdinov, A. M. Semikhatov and I. Yu. Tipunin, 
{ Modular group representations and fusion in logarithmic conformal field theories and in the quantum group center,} 
{\it Comm. Math. Phys}. {\bf 265} (2006)  47--93. 
%\bibitem[FGST2]{F2}
%Feigin, B.L.;  Gainutdinov, a. M.; Semikhatov, A. M.;  Tipunin, I. Yu.
%{\it Kazhdan--Lusztig correspondence for the representation category of the triplet W-algebra in logarithmic CFT,}
%math.QA/0512621.  
\bibitem{H}
G. Hennings,  
{ Invariants of links and 3-manifods obtained from Hopf algebras},
{\it J. London Math. Soc}. (2) {\bf 54} (1996) 594--624.  
\bibitem{J}
M. Jimbo, T. Miwa and Y.  Takeyama, 
{Counting minimal form factors of the restricted sine-Gordon model,}
{\it Mosc. Math. J}. {\bf 4} (2004) 787--846, 981.
\bibitem{Ker}
T. Kerler, 
{Mapping class group actions on quantum doubles},
{\it Comm. Math. Phys}. {\bf 168} (1995) 353--388.  
%\bibitem{K}
%Kassel, C.
%{Quantum Groups, }
%1995, Springer-Verlag, New York.  
%\bibitem{MY}
%Murakami, H. and Yokota, Y.
%{\it The colored Jones polynomials of the figure-eight knot and its Dehn surgery spaces}, 
%J. Reine Angew. Math. {\bf 607} (2007), 47--68.
\bibitem{KR}
L. Kauffman and D. Radford,   
{Invariants of 3-manifods derived from finite-dimensional Hopf algebras},
{\it J. Knot Theory Ramifications} {\bf 4} (1995) 131--162.  
\bibitem{Ku}
G. Kuperberg, 
{Involutory Hopf lagebras and 3-manifod invariants},
{\it Internat. J. Math}. {\bf 2} (1991) 41--66.  
\bibitem{M}
J. Murakami, 
{Colored Alexander invariants and cone manifolds}, 
to appear in {\it Osaka J. Math}. 
\bibitem{O}
T. Ohtsuki, 
{Colored ribbon Hopf algebas and universal invariants of framed links,}
{\it J. Knot Theory Ramifications} {\bf 2} (1993)  211--232.
\bibitem{O2}
T. Ohtsuki, 
{Invariants of 3-manifods derived from universal invariants of framed links},
{\it Math. Proc. Cambridge Phillos, Soc}.  {\bf 117} (1995)  259--273.
\bibitem{T}
V. Turaev, 
{\it Quantum invariants of knots and 3-manifolds},
(de Gruyter, Berlin, 1994).  
\end{thebibliography}
\end{document}